\documentclass[twoside,11pt,headrule]{amsart}
\usepackage{enumerate, latexsym, amsmath, amsfonts, amssymb, amsthm, color}
\def\pmod #1{\ ({\rm{mod}}\ #1)}
\def\Z{\Bbb Z}
\def\N{\Bbb N}

\def\l{\left}
\def\r{\right}
\def\bg{\bigg}
\def\({\bg(}
\def\){\bg)}
\def\t{\text}
\def\f{\frac}

\def\ls{\leqslant}
\def\gs{\geqslant}

\def\sm{\setminus}
\def\bi{\binom}
\def\al{\alpha}

\def\ve{\varepsilon}

\def\eq{\equiv}

\def\nh{\noalign{\hrule}}
\def\hh{height4pt}

\def\om{&\omit &}
 
 \def\Proof{\noindent{\it Proof}}

 \def\Ack{\medskip\noindent {\bf Acknowledgments}}
\def\FF#1#2#3{{}_2F_1\bigg(\bmatrix{#1}\\{#2}\end{bmatrix}\bigg|#3\bigg)}

\theoremstyle{plain}
\newtheorem{theorem}{Theorem}

\newtheorem{corollary}{Corollary}
\newtheorem{conjecture}{Conjecture}
\theoremstyle{definition}

\theoremstyle{remark}
\newtheorem{remark}{Remark}
\renewcommand{\theequation}{\arabic{section}.\arabic{equation}}
\renewcommand{\thetheorem}{\arabic{section}.\arabic{theorem}}
\renewcommand{\thelemma}{\arabic{section}.\arabic{lemma}}
\renewcommand{\thecorollary}{\arabic{section}.\arabic{corollary}}
\renewcommand{\theconjecture}{\arabic{conjecture}}
\renewcommand{\theremark}{\arabic{remark}}

\vspace{4mm}
\begin{document}
\renewcommand{\baselinestretch}{1.3}
\renewcommand{\arraystretch}{1.3}

\hbox{Nanjing Univ. J. Math. Biquarterly 36 (2019), no.\,2, 108-133.}
\medskip

\title
[{New observations on primitive roots modulo primes}]
{New observations on primitive roots modulo primes}

\author
[Zhi-Wei Sun] {Zhi-Wei Sun}
\address {Department of Mathematics, Nanjing
University, Nanjing 210093, People's Republic of China}
\email{zwsun@nju.edu.cn}

\thanks{
\newline \indent Supported by the National Natural Science Foundation (grant 11571162)
 of China.}

 \begin{abstract} We make many new observations on primitive roots modulo primes. For an odd prime $p$ and an integer $c$, we establish a theorem concerning
$\sum_g(\frac{g+c}p)$, where $g$ runs over all the primitive roots modulo $p$ among $1,\ldots,p-1$,
and $(\frac{\cdot}p)$ denotes the Legendre symbol.
On the basis of our numerical computations, we formulate 35 conjectures
involving primitive roots modulo primes. For example, we conjecture that for any prime $p$
  there is a primitive root $g<p$ modulo $p$ with $g-1$ a square, and that for any prime $p>3$ there is a prime $q<p$ with
 the Bernoulli number $B_{q-1}$ a primitive root modulo $p$. We also make related observations on
quadratic nonresidues modulo primes and  primitive prime divisors of some combinatorial sequences. For example, based on heuristic arguments we conjecture that
 for any prime $p>3$ there exists a Fibonacci number $F_k<p/2$ which is a quadratic nonresidue modulo $p$; this implies that there is a deterministic polynomial time algorithm to find square roots of quadratic residues modulo a prime $p>3$.
\end{abstract}

\maketitle

\renewcommand{\baselinestretch}{1.3}
\renewcommand{\arraystretch}{1.3}

\vskip 0.5cm
%\section{Introduction}
 \noindent{\Large\bf 1\quad Introduction}
 \vskip 0.5cm

\renewcommand{\theequation}{1.\arabic{equation}}
\renewcommand{\thetheorem}{1.\arabic{theorem}}
\renewcommand{\thecorollary}{1.\arabic{corollary}}
\renewcommand{\thelemma}{1.\arabic{lemma}}
\renewcommand{\theremark}{1.\arabic{remark}}
\setcounter{lemma}{0} \setcounter{theorem}{0}
\setcounter{corollary}{0} \setcounter{remark}{0}
\setcounter{equation}{0} \setcounter{conjecture}{0}

 Let $p$ be any prime. It is well known that $\mathbb F_p=\Z/p\Z=\{\bar a=a+p\Z:\ a\in\Z\}$ is a field and $\mathbb F_p^*=\mathbb F_p\sm\{0\}$
 is a cyclic group of order $p-1$. A rational $p$-adic integer $g$ is called a {\it primitive root} modulo $p$ if $\bar g =g\ \t{mod}\ p$
 is a generator of $\mathbb F_p^*$. The standard proof of the existence of a primitive root modulo $p$ (cf. [IR, p.\,40]) is nonconstructive, and it provides no way to find
 an explicit primitive root modulo $p$.

 The most famous unsolved problem on primitive roots modulo primes is the following conjecture posed by E. Artin in 1927 (see [M] for a survey of results
 towards Artin's conjecture).

\medskip
 \noindent{\bf Artin's Conjecture}\quad  {\it If $g\in\Z$ is neither $-1$ nor a square, then there are infinitely many primes $p$ such that $g$ is a primitive root modulo $p$.}
 \medskip

 Let $p$ be a prime. It is well known that the set
 \begin{equation}\label{G(p)}G(p):=\{g\in\{1,\ldots,p-1\}:\ g \ \t{is a primitive root modulo}\ p\}|\end{equation}
 has cardinality $\varphi(p-1)$, where $\varphi$ denotes Euler's totient function. According to \cite[p.\,377]{Gu}, P. Erd\"os ever asked the following open question.

 \medskip
 \noindent{\bf Erd\H os' Problem}\quad  {\it  Whether for any sufficiently large prime $p$ there exists a prime $q<p$ which is a primitive root modulo $p$?}
 \medskip

 Let $q>1$ be a prime power. For the finite field $\mathbb F_q$ of order $q$, the multiplicative group $\mathbb F_q^*=\mathbb F_q\sm\{0\}$ is a cyclic group of order $q-1$
 and any generator of this group is called a primitive root (or primitive element) of the field $\mathbb F_q$.
 In 1971 E. Vegh \cite{V} guessed that if $q>61$ then any element of $\mathbb F_q$ can be written as a difference of two primitive roots of $\mathbb F_q$.
 In 1984 S. W. Golomb \cite{G} conjectured that any nonzero element of $\mathbb F_q$ can be expressed as a sum of two primitive roots of $\mathbb F_q$.
 After many earlier efforts to prove Vegh's and Golomb's conjectures and their linear extensions, it is now known that
 if $q>61$ and $a,b,c\in \mathbb F_q^*$ then there always exist primitive roots $g$ and $h$ of $\mathbb F_q$ with $a=bg+ch$ (see the introduction part of the recent paper \cite{COT}).
In particular, this implies that for any prime $p>61$ the set $G(p)$ defined in (1.1) contains two consecutive integers.
In contrast, the twin prime conjecture still remains unsolved despite Y. Zhang's breakthrough (cf. \cite{Z}) on prime gaps.

In 1989 W. B. Han \cite{H} studied extensions of Vegh's and Golomb's conjectures to polynomials over finite fields. Using Weil's theorem on character sums,
he established the following general theorem.

\begin{theorem}\label{Th1.1} {\rm (Han [H])} Let $q>1$ be a prime power. Let $f(x)$ and $g(x)$ be polynomials over the finite field $\mathbb F_q$ such that none of $g(x)$
and $f(x)g(x)^k\ (k=0,1,2,\ldots)$ can be written in the form $cP(x)^d$ with $c\in\mathbb F_q$, $1<d\mid(q-1)$ and $P(x)\in \mathbb F_q[x]$.
Let $m$ be the number of distinct zeroes of $f(x)$ in the splitting field of $f(x)$, and let $n$ be the number of distinct zeroes
of $g(x)$ in the splitting field of $g(x)$. If $\sqrt q\gs (m+n-1)4^{\omega(q-1)}$, then for some $a\in \mathbb F_q$ both $f(a)$ and $g(a)$ are primitive roots of $\mathbb F_q$,
where $\omega(q-1)$ the number of distinct prime divisors of $q-1$.
\end{theorem}

As a consequence of Theorem \ref{Th1.1}, Han noted that for any finite field $\mathbb F_q$ with $q\gs2^{66}$, if $a,b,c\in\mathbb F_q$ and $ac(b^2-4ac)\not=0$,
then there is a primitive root $g\in \mathbb F_q$ with $ag^2+bg+c$ also a primitive root of $\mathbb F_q$ (cf. \cite[Corollary 3]{H}).
In particular, for any prime $p>2^{66}$ there is a primitive root $g$ modulo $p$ such that $g^2+1$ is also a primitive root modulo $p$.
In contrast, it is unproven that there are infinitely many primes of the form $x^2+1$ with $x\in\Z$.

On Oct. 3, 2013, the author \cite[A229910]{S} conjectured that
for any prime $p>13$ there is a primitive root $g$ modulo $p$
such that $g+g^{-1}$ is also a primitive root modulo $p$.
Based on Han's work, the author  showed in Oct. 2013 (cf. \cite[A229910]{S}) that for each $\ve\in\{\pm1\}$ and
for any finite field $\mathbb F_q$ with $q>2^{66}$, there is a primitive element $g$ of $\mathbb F_q$ such that $g + \ve g^{-1}$ is also a primitive root of $\mathbb F_q$.

In 2018, S.D. Cohen, T. Oliveira e Silva and Sutherland \cite{COS} obtained the following further result.

\begin{theorem}\label{Th1.2} {\rm (\cite[Corollary 2]{COS})} Let $q>5$ be a prime power.

{\rm (i)} If $q\not\in\{7,9,13,25,121\}$, then there is a primitive element $g$ of $\mathbb F_q$
with $g+g^{-1}$ also primitive.

{\rm (ii)} If $q\not\in\{9,13,25,61,121\}$, then there is a primitive element $g$ of $\mathbb F_q$
with $g-g^{-1}$ also primitive.
\end{theorem}

For any odd prime $p$ and integer $c$, we introduce
\begin{equation}\label{Def-S}S_p(c):=\sum_{g\in G(p)}\l(\f{g+c}p\r),
\end{equation}
where $(\f{\cdot}p)$ denotes the Legendre symbol.
Note that if $p\mid c$ then $S_p(c)=-|G(p)|=-\varphi(p-1)$.
Concerning $S_p(c)$ we have the following result.

\begin{theorem}\label{Th-g} Let $p$ be any odd prime.

{\rm (i)} We have
\begin{equation}\label{Sp1}S_p(1)=0.
\end{equation}

{\rm (ii)} For any integer $c\in\Z$ with $c\not\eq0\pmod p$, we have
\begin{equation}\label{Spc} S_p(c)\eq\l(\f cp\r)\sum_{k=0}^{(p-1)/2}\f{\bi{2k}k}{(-4c)^k}\mu\l(\f{p-1}{(k,p-1)}\r)\f{\varphi(p-1)}{\varphi((p-1)/(k,p-1))}\pmod p,\end{equation}
where $\mu$ is the M\"obius function
and $(k,p-1)$ is the greatest common divisor of $k$ and $p-1$.
\end{theorem}
\Proof. (i) For $g\in\{1,\ldots,p-1\}$ let $g^*\in\{1,\ldots,p-1\}$
be the inverse of $g$ modulo $p$ (i.e., $gg^*\eq1\pmod p$). Clearly,
\begin{align*}S_p(1)=&\sum_{g\in G(p)}\l(\f{g^*+1}p\r)=-\sum_{g\in G(p)}\l(\f{g(g^*+1)}p\r)
=-\sum_{g\in G(p)}\l(\f{1+g}p\r)=-S_p(1)
\end{align*}
and thus $S_p(1)=0$.

(ii) Now let $c\in\Z$ with $c\not\eq0\pmod p$. Observe that
\begin{align*} \sum_{g\in G(p)}(g+c)^{(p-1)/2}=&\sum_{g\in G(p)}\sum_{k=0}^{(p-1)/2}\bi{(p-1)/2}kg^k c^{(p-1)/2-k}
\\=&\sum_{k=0}^{(p-1)/2}\bi{(p-1)/2}kc^{(p-1)/2-k}\sum_{g\in G(p)}g^k
\\\eq&\l(\f cp\r)\sum_{k=0}^{(p-1)/2}\f{\bi{-1/2}k}{c^k}\sum_{g\in G(p)}g^k
\\=&\l(\f cp\r)\sum_{k=0}^{(p-1)/2}\f{\bi{2k}k}{(-4c)^k}\sum_{g\in G(p)}g^k\pmod p.
\end{align*}
Fix a primitive root $g_0$ modulo $p$. Then
$$\sum_{g\in G(p)}g^k\eq\sum^{p-1}_{j=1\atop (j,p-1)=1} g_0^{jk}\eq\varphi(p-1)\f{\mu((p-1)/(k,p-1))}{\varphi((p-1)/(k,p-1))}\pmod p$$
via the known evaluations of Ramanujan sums.
Therefore the desired \eqref{Spc} follows.

In view of the above, we have completed the proof of Theorem \ref{Th-g}. \qed

\begin{corollary}\label{Cor-g} {\rm (i)} For any prime $p\eq1\pmod 4$, we have
$S_p(-c)=S_p(c)$ for all $c\in\Z$, in particular $S_p(-1)=S_p(1)=0.$

{\rm (ii)} Let $p$ be a Fermat prime and let $c\in\Z$ with $c\not\eq0\pmod p$. Then
\begin{equation}\label{Fermat}S_p(c)=\f{1-(\f cp)}2.
\end{equation}
\end{corollary}
\Proof. (i) We now prove the first part. Let $c\in\Z$. If $p\mid c$, then $S_p(-c)=S_p(0)=S_p(c)$.

Now we assume $p\nmid c$.
If $k$ is odd, then $(p-1)/(k,p-1)\eq0\pmod 4$ and hence $\mu((p-1)/(k,p-1))=0$. If $k$ is even, then
$(-c)^k=c^k$. So, by \eqref{Spc} we have $S_p(-c)\eq S_p(c)\pmod p$. Since
$$|S_p(\pm c)|\ls|G(p)|=\varphi(p-1)\ls\f{p-1}2,$$
we must have $S_p(-c)=S_p(c)$. In particular, $S_p(-1)=S_p(1)=0$ with the aid of \eqref{Sp1}.

(ii) Now we turn to show the second part of Corollary \ref{Cor-g}. Write $p=2^n+1$ with $n$ a power of two. For any positive integer $k<(p-1)/2=2^{n-1}$, we have
$(k,p-1)\mid 2^{n-2}$ and hence $\mu((p-1)/(k,p-1))=0$ since $(p-1)/(k,p-1)$ is divisible by $4$.
Thus, by \eqref{Spc} we have
\begin{align*}\l(\f cp\r)S_p(c)\eq&\varphi(p-1)+\f{\bi{p-1}{(p-1)/2}}{(-4c)^{(p-1)/2}}\mu(2)\f{\varphi(p-1)}{\varphi(2)}
\\\eq&\varphi(p-1)\l(1-\l(\f cp\r)\r)=\f{p-1}2\l(1-\l(\f cp\r)\r)\eq\f{(\f cp)-1}2\pmod p.
\end{align*}
Note that $|S_p(c)|\ls \varphi(p-1)=(p-1)/2$. Therefore \eqref{Fermat} holds.
\qed

In view of Erd\H os' problem, Theorems 1.1-1.2 and various problems on primes of special forms,
we are led to consider whether primitive roots modulo primes can take certain special forms.
In Section 3 we will pose various conjectures in this direction based on our computational checks.
Since any primitive root modulo an odd prime $p$ must be a quadratic nonresidue modulo $p$,
in Section 2 we will investigate quadratic nonresidues (modulo primes) of certain special forms
armed with heuristic arguments. In Section 4, we will pose some other conjectures involving primitive roots modulo primes.

Let $(a_n)_{n\gs1}$ be a sequence of integers. If no term of the sequence $(a_n)_{n\gs1}$ has a prime divisor greater than a given integer $N>1$, then for any prime $q\eq1\pmod{4\prod_{p\ls N}p}$
we have $(\f{a_n}q)=1$ for all $n=1,2,3,\ldots$. If a prime $p$ divides the $n$-th term $a_n$ but it does not divide any previous term $a_k$ with $0<k<n$, then $p$ is called a {\it primitive prime divisor} of the term $a_n$. For our purposes, we are interested in those integer sequences
with infinitely many terms having primitive prime divisors.

In 1886 A. S. Bang \cite{B} proved that for any integer $n>1$ with $n\not=6$ the number $2^n-1$ has a prime divisor not dividing any $2^k-1$ with $k\in\{1,\ldots,n-1\}$.
In 1892 K. Zsigmondy \cite{Zs} extended this as follows: If $a$ and $b$ are integers with $a>b>0$ and $(a,b)=1$, then for any integer $n>2$ the number $a^n-b^n$ has a prime divisor not dividing any $a^k-b^k$ with $0<k<n$, unless $a=2$, $b=1$ and $n=6$.

Recall that the Fibonacci numbers are given by
$$F_0=0,\ F_1=1,\ \t{and}\ F_{n+1}=F_n+F_{n-1}\ (n=1,2,3,\ldots).$$
Carmichael's theorem (cf. \cite{C}) asserts that for any integer $n>12$ the $n$-th Fibonacci number $F_n$
has a prime divisor $p$ which does not divide any previous Fibonacci number $F_k$ with $0<k<n$.
 Let $A,B\in\Z$ with $B\not=0$ and $A^2\not=4B$. The Lucas sequence $u_n=u_n(A,B)\ (n=0,1,2,\ldots)$ is defined by
$$u_0=0,\ u_1=1,\ \t{and}\ u_{n+1}=Au_n-Bu_{n-1}\ \t{for}\ n=1,2,3,\ldots.$$
In 2001 Y. Bilu, G. Hanrot and P. M. Voutier \cite{BHV} proved that for any integer $n>30$ the term $u_n(A,B)$ has prime divisor not dividing any previous term $u_k(A,B)$ with $0<k<n$.

In Section 5 we look at various combinatorial sequences of integers or rational numbers to see whether larger terms have primitive prime divisors.
This leads us to generate some tables on primitive prime divisors and formulate various conjectures in this direction.

Throughout this paper, we set $\N=\{0,1,2,\ldots\}$ and $\Z^+=\{1,2,3,\ldots\}$.

\vskip 0.5cm
 \noindent{\Large\bf 2\quad On Special Quadratic Nonresidues Modulo Primes}
 \vskip 0.5cm

\renewcommand{\theequation}{2.\arabic{equation}}
\renewcommand{\thetheorem}{2.\arabic{theorem}}
\renewcommand{\thecorollary}{2.\arabic{corollary}}
\renewcommand{\thelemma}{2.\arabic{lemma}}
\renewcommand{\theconjecture}{2.\arabic{conjecture}}

 \setcounter{lemma}{0} \setcounter{theorem}{0}
\setcounter{corollary}{0}
\setcounter{remark}{0}
\setcounter{equation}{0}
\setcounter{conjecture}{0}

Let $a$ be a quadratic residue modulo an odd prime $p$. How to solve the congruence $x^2\eq a\pmod p$ efficiently?
By the Tonelli-Shanks Algorithm (cf. R. Crandall and C. Pomerance \cite[pp. 93-95]{CP}), if we know a quadratic nonresidue $d\in\Z$ modulo $p$ then
one can solve $x^2\eq a\pmod p$ efficiently as follows:

Write $p-1=2^st$ with $s,t\in\Z^+$ and $2\nmid t$, and find even integers $m_1,\ldots,m_s$ with $(ad^{m_i})^{2^{s-i}t}\eq1\pmod p$
for all $i=1,\ldots,s$ in the following way:
$m_1:=0$, and after those $m_1,\ldots, m_i$ (with $1\ls i<s$) have been chosen we select $m_{i+1}\in\{m_i,m_i+2^i\}$ such that $(ad^{m_{i+1}})^{2^{s-i-1}t}\eq1\pmod p$.
Note that $((ad^{m_i})^{2^{s-i-1}t})^2\eq1\pmod p$ and hence $(ad^{m_i})^{2^{s-i-1}t}\eq\pm1\pmod p$. If $(ad^{m_i})^{2^{s-i-1}t}\eq-1\pmod p$, then
$$(ad^{m_i+2^i})^{2^{s-i-1}t}\eq-d^{2^{s-1}t}=-d^{(p-1)/2}\eq1\pmod p.$$
As $(ad^{m_s})^t\eq1\pmod p$, we have $x^2\eq a\pmod p$ with $x=\pm a^{(t+1)/2}(d^t)^{m_s/2}$.

However, there is no known deterministic, polynomial time algorithm for finding a quadratic nonresidue $d$
modulo a given odd prime $p$. According to \cite[pp.\,93-95]{CP}, under the Extended Riemann Hypothesis for algebraic fields, it can be shown that there is a positive quadratic nonresidue $d<2\log^2p$; and so an exhaustive search to this limit
succeeds in finding a quadratic nonresidue in polynomial time. Thus, under the ERH, one can find square roots for quadratic residues modulo the prime $p$ in deterministic, polynomial time.

As the Fibonacci numbers grow exponentially, part (i) of our following conjecture is particularly interesting since
it implies that we can find square roots for quadratic residues modulo a prime $p>3$ in deterministic polynomial time.

\begin{conjecture}\label{Conj2.1} {\rm (i) (2014-04-26)} For any integer $n>4$, there is a Fibonacci number $f<n/2$ with $x^2\eq f\pmod n$ for no integer $x$.

{\rm (ii) (2014-04-27)} For any odd prime $p$, let $f(p)$ be the least Fibonacci number with $(\f{f(p)}p)=-1$. Then $f(p)=o(p^{0.7})$ as $p\to\infty$.
Moreover, we have $f(p)=O(p^c)$ for any $c>c_0=\log_2\f{1+\sqrt5}2\approx 0.694$.

{\rm (iii) (2014-05-07)} For any prime $p$, there exists a positive integer $k\ls\sqrt{p+2}+2$ such that $F_k+1$ is a primitive root modulo $p$.
\end{conjecture}

Conjecture 2.1(i) can be reduced to the case when $n$ is prime. In fact, for any positive integer $n$ divisible by $3$ or $4$, there is no square congruent to $F_3=2$ modulo $n$.
If $n>4$ has a prime divisor $p>3$, and there is a positive Fibonacci number $F_k<p/2$ with $x^2\not\eq F_k\pmod p$ for all $x\in\Z$, then
$F_k<n/2$ and also $x^2\not\eq F_k\pmod n$ for all $x\in\Z$. We have verified part (i) for every $n=4,5,\ldots,3\times10^9$.
For data and graphs related to Conjecture 2.1(i), one may consult \cite[A241568, A241604 and A241675]{S}.
\medskip

As for part (ii) of Conjecture 2.1, we don't have a rigorous proof but it seems reasonable in view of the following heuristic arguments.

\medskip
\noindent {\bf Heuristic Arguments for Conjecture 2.1(ii)}.
In light of Carmichael's theorem on primitive prime divisors of Fibonacci numbers, we may think that a positive Fibonacci number not exceeding $p^c$
is a quadratic residue modulo $p$ with `probability' $1/2$. Roughly speaking, there are about
$$\f{\log_2p^c}{\log_2\f{1+\sqrt5}2} = \f c{c_0}\log_2p$$
positive Fibonacci numbers not exceeding $p^c$. So we might expect that all positive Fibonacci numbers not exceeding $p^c$
are quadratic residues modulo $p$ with probability
$$\l(\f12\r)^{(\log_2p)c/c_0}=\f1{p^{c/c_0}}.$$
As $\sum_p p^{-c/c_0}$ converges, it seems reasonable to think that there are finitely many primes $p$ for which all positive Fibonacci numbers
not exceeding $p^c$ are quadratic residues modulo $p$. So the guess $f(p)=O(p^c)$ probably holds.

\medskip

We have verified Conjecture 2.1(iii) for all primes $p<10^8$, and observed that no Fibonacci number is a primitive root modulo the prime $3001$.
Note that for any integer $n>1$ there is a Fibonacci number $F_k$ with $F_k+1\eq0\pmod n$. In fact, by the Pigeonhole Principle,
there are $0\ls i<j\ls n^2$ such that $F_i\eq F_j\pmod n$ and $F_{i+1}\eq F_{j+1}\pmod n$, and hence $F_{j-i}\eq F_0=0\pmod n$ and $F_{j-i+1}\eq F_1=1\pmod n$.
Clearly $k=j-i-2>0$ since $F_1=F_2=1\not\eq0\pmod n$, and
$$F_k=F_{k+2}-(F_{k+3}-F_{k+2})=2F_{k+2}-F_{k+3}=2F_{j-i}-F_{j-i+1}\eq-1\pmod n.$$
\medskip

Recall that the Lucas numbers $L_0,L_1,L_2,\ldots$ are defined by
$$L_0=2,\ L_1=1,\ \t{and}\ L_{n+1}=L_n+L_{n-1}\ (n=1,2,3,\ldots).$$
It is well known that
$$L_n=\l(\f{1+\sqrt5}2\r)^n+\l(\f{1-\sqrt5}2\r)^n\quad\t{for all}\ n\in\N.$$
Our following conjecture is similar to Conjecture 2.1.

\begin{conjecture}\label{Conj2.2} {\rm (i) (2014-04-26)} For any integer $n>2$, there is a Lucas number $L_k<n$ such that $x^2\not\eq L_k\pmod n$ for all $x\in\Z$.

{\rm (ii) (2014-04-27)}  For any odd prime $p$, let $\ell(p)$ be the least Lucas number with $(\f{\ell(p)}p)=-1$. Then $\ell(p)=o(p^{0.7})$ as $p\to\infty$.
Moreover, we have $\ell (p)=O(p^c)$ for any $c>\log_2\f{1+\sqrt5}2\approx 0.694$.

{\rm (iii) (2014-05-21)} For any prime $p$, there exists a positive integer $k<\sqrt{p}+2$ such that $L_k+1$ is a primitive root modulo $p$.
\end{conjecture}

We have verified Conjecture \ref{Conj2.2}(i) for every $n=3,\ldots,10^9$, and Conjecture \ref{Conj2.2}(iii) for all primes $p<10^8$. The least Lucas number which is a quadratic nonresidue
modulo the prime $p=167$, is $L_{10}=123>167/2$. Note that no Lucas number is a primitive root modulo the prime $28657$. Also, for any integer $n>1$ there is a positive integer $j<n^2$ such that $L_j\eq L_0=2\pmod n$
and $L_{j+1}\eq L_1=1\pmod n$, and hence $L_{j-1}=L_{j+1}-L_j\eq 1-2=-1\pmod n$.

The following conjecture similar to Conjectures 2.1 and 2.2 is concerned with cubic nonresidues modulo primes.
For a prime $p\eq1\pmod 3$, it seems reasonable to think that $2^k-1$ is a cubic nonresidue modulo $p$ with probability $2/3=1/1.5$.

\begin{conjecture}\label{Conj2.3} {\rm (2014-05-11)} Let $p$ be any prime with $p\eq1\pmod3$. Then, there is a positive integer $k$ with $2^k-1<p/2$ such that $2^k-1$ is a cubic nonresidue modulo $p$.
Moreover, for any $c>\log1.5/\log2\approx 0.585$ we have $s(p)=O(p^c)$, where $s(p)$ denotes the least positive cubic nonresidue modulo $p$ in the form $2^k-1$ with $k\in\Z^+$.
\end{conjecture}

We have verified the first assertion in Conjecture 2.3 for all primes $p<10^9$ with $p\eq1\pmod3$; for example,
the least positive cubic nonresidue modulo the prime $p=4667629$ in the form $2^k-1$ is $2^{15}-1=32767$.
The second assertion in Conjecture 2.3 sounds reasonable by heuristic arguments.

To conclude this section, we pose one more conjecture.

\begin{conjecture}\label{Conj2.4} {\rm (2014-04-20)}
For any prime $p>7$, there is a prime $q<p$ with $2^q-1$ a quadratic residue modulo $p$.
Also, for each prime $p>5$, there exists a prime $q<p$ such that $2^q+1$ is a quadratic nonresidue modulo $p$.
\end{conjecture}

We have verified Conjecture 2.4 for primes $p$ below $10^8$;
see \cite[A235709 and A235712]{S} for related data and graphs.
For example, $41$ is  the least prime $q<13003$ with $2^q-1$ a quadratic residue modulo the prime $13003$.
Note that for the prime $p=2089$ there is no prime $q<p$ with $2^q+1$ a primitive root modulo $p$.

\vskip 0.5cm
 \noindent{\Large\bf 3\quad On Primitive Roots of Special Forms}
 \vskip 0.5cm

\renewcommand{\theequation}{3.\arabic{equation}}
\renewcommand{\thetheorem}{3.\arabic{theorem}}
\renewcommand{\thecorollary}{3.\arabic{corollary}}
\renewcommand{\thelemma}{3.\arabic{lemma}}
\renewcommand{\theconjecture}{3.\arabic{conjecture}}
\renewcommand{\theremark}{3.\arabic{remark}}

 \setcounter{lemma}{0} \setcounter{theorem}{0}
\setcounter{corollary}{0}
\setcounter{remark}{0}
\setcounter{equation}{0}
\setcounter{conjecture}{0}

As we mentioned in Section 1, it is known that for any sufficiently large prime $p$ there is a primitive root modulo $p$
in the form $x^2+1$ with $x\in\Z$. Part (i) of our following conjecture is stronger than this.

\begin{conjecture}\label{Conj3.1} {\rm (2014-04-23) (i)} Every prime $p$ has a primitive root $g
< p$ modulo $p$ of the form $k^2 + 1$. In other words, for any prime $p$, there is a
primitive root $0 < g < p$ modulo $p$ with $g - 1$ an integer square.

{\rm (ii)}  For any prime $p > 3$, there is a triangular number $g < p$ which is a
primitive root modulo $p$. Also, every prime $p > 7$ has a primitive root
 $g<p$ modulo $p$ which is a product of two consecutive integers.
\end{conjecture}

\begin{remark}\label{Rem3.1}\rm The author verified Conjecture 3.1(i) for all primes
$p<10^7$ in April 2014, and later C. Greathouse \cite{Gr} extended the verification
to all primes below $p<10^{10}$ in May 2014. The author would like to offer 2,000 RMB
as the prize for the first complete solution of Conjecture 3.1(i).
Note that for any prime $p>5$, one of the three numbers $1^2+1=2$,
$2^2+1=5$ and $3^2+1=10=2\times5$ is a quadratic residue modulo $p$.
We have verified Conjecture 3.1(ii) for primes $p<10^9$. Note that for any prime $p>3$ one of
$1\times2,\ 2\times3$ and $3\times 4$ is a quadratic residue modulo $p$.
For data and graphs concerning Conjecture 3.1, one may consult \cite[A239957, A241476, A239963 and A241492]{S}.
\end{remark}

\medskip
\centerline{Table 3.1: Primes $p$ with unique primitive root $g$ of the form $k^2+1<p$}
\smallskip
\centerline{\vbox{\offinterlineskip
\halign{\vrule#&\ \ #\ \hfill   &&\vrule#&\ \ \hfill#\ \
\cr\nh\cr \hh\om\om\om\om\om\om\om\om\om\om\om
\cr &$p$& &$2$& &$3$& &$5$& &$7$& &$11$&  &$13$& &$31$& &$71$& &$79$& &$151$&
\cr \hh\om\om\om\om\om\om\om\om\om\om\om
\cr\nh\cr \hh\om\om\om\om\om\om\om\om\om\om\om
\cr &$k$& &$1$& &$1$& &$1$& &$2$& &$1$& &$1$& &$4$& &$8$& &$6$& &$9$&
\cr \hh\om\om\om\om\om\om\om\om\om\om\om
\cr\nh\cr \hh\om\om\om\om\om\om\om\om\om\om\om
\cr &$g=k^2+1$& &$2$& &$2$& &$2$& &$5$& &$2$& &$2$& &$17$& &$65$& &$37$& &$82$&
\cr \hh\om\om\om\om\om\om\om\om\om\om\om
\cr \nh\cr }}}
\smallskip
\bigskip

In 2000 D.K.L. Shiu \cite{Sh} proved that if $a\in\Z$ and $m\in\Z^+$ are relatively prime then for any $k\in\Z^+$
there is a positive integer $n$ such that $p_{n+1}\eq p_{n+2}\eq\ldots\eq p_{n+k}\eq a\pmod m$, where $p_j$ denotes the $j$-th prime. This remarkable result
implies that the set $\{S_n=\sum_{k=1}^n p_k:\ n=1,2,3,\ldots\}$ contains a complete system of residues modulo any positive integer $m$.
In \cite{S13a} the author conjectured that the set $\{s_n=\sum_{k=1}^n (-1)^{n-k}p_k:\ n=1,2,3,\ldots\}$ also contains a complete system of residues modulo any positive integer $m$.
Motivated by these, we pose the following conjecture.

\begin{conjecture}\label{Conj3.2} {\rm (2014-05-10) (i)} For any odd prime $p$, there is a primitive root $g<p$ modulo $p$ in the form $S_n=\sum_{k=1}^np_k$ with $n\in\Z^+$.

{\rm (ii)} For any integer $n>1$, there is a number $k\in\{1,\ldots,n\}$ such that $s_k=\sum_{j=1}^k(-1)^{k-j}p_j$ is a primitive root modulo $p_n$.

{\rm (iii)} For any integer $n>4$, there is a positive integer $k\ls n/2$ such that
$\prod_{j=1}^kp_j$ is a primitive root modulo $p_n$.
\end{conjecture}
\begin{remark}\label{Rem3.2} We have verified part (i) for all odd primes $p<10^9$, and parts (ii) and (iii) for $n$ up to $10^7$.
See \cite[A242266 and A242277]{S} for related data and graphs.
\end{remark}

\newpage
\medskip
\centerline{Table 3.2: Primes $p$ with unique primitive root $g$ of the form $\sum_{k=1}^np_k<p$}
\smallskip
\centerline{\vbox{\offinterlineskip
\halign{\vrule#&\ \ #\ \hfill   &&\vrule#&\ \ \hfill#\ \
\cr\nh\cr \hh\om\om\om\om\om\om\om\om\om\om
\cr &$p$& &$3$& &$5$& &$7$& &$11$& &$13$&  &$31$& &$71$& &$127$& &$241$&
\cr \hh\om\om\om\om\om\om\om\om\om\om
\cr\nh\cr \hh\om\om\om\om\om\om\om\om\om\om
\cr &$n$& &$1$& &$1$& &$2$& &$1$& &$1$& &$4$& &$5$& &$7$& &$10$&
\cr \hh\om\om\om\om\om\om\om\om\om\om
\cr\nh\cr \hh\om\om\om\om\om\om\om\om\om\om
\cr &$g=\sum_{k=1}^np_k$& &$2$& &$2$& &$5$& &$2$& &$2$& &$17$& &$28$& &$58$& &$129$&
\cr \hh\om\om\om\om\om\om\om\om\om\om
\cr \nh\cr }}}
\smallskip
\bigskip

For each $n\in\Z^+$ let $p(n)$ be the number of ways to write $n$ as a sum of some unordered positive integers with repetitions allowed. This is the well-known partition function.
On April 24, 2014 the author conjectured that for any prime $q$
there is a positive integer $n$ with $p(n)<q$ such that $p(n)$ is a primitive root modulo $q$
(cf. \cite[Conjecture 4.10(i)]{S14}).
We have verified this for all primes $q<10^9$.

\begin{conjecture}\label{Conj3.3} {\rm (i) (2014-04-22)} For any prime $p > 3$, there exists a prime $q < p/2$ such that the
Mersenne number $M_q = 2^q - 1$ is a primitive root modulo $p$.

{\rm (ii) (2014-05-09)} For any prime $p>3$, there exists a positive integer $g<p$
such that $g$, $2^g-1$ and $(g-1)!$ are all primitive roots modulo $p$.
\end{conjecture}

\begin{remark}\label{Rem3.3}\rm (a) We have verified Conjecture 3.3(i) for all primes
$3<p < 10^7$; see \cite[A236966]{S} for related data and graphs. For example, for the prime
$p = 5336101$, the least prime $q < p/2$ with $2^q - 1$ a primitive root
modulo $p$ is 193.

(b) Conjecture 3.3(ii) is very strong! We have verified it for all primes $3<p<10^7$;
see \cite[A242248 and A242250]{S}
for related data and graphs.
\end{remark}

\centerline{Table 3.3: Primes $p$ with unique $0<g<p$ such that}
\centerline{$g$, $2^g-1$ and $(g-1)!$ are all primitive roots mod $p$}
\smallskip
\centerline{\vbox{\offinterlineskip
\halign{\vrule#&\ \ #\ \hfill   &&\vrule#&\ \ \hfill#\ \
\cr\nh\cr \hh\om\om\om\om\om\om\om\om\om\om\om
\cr &$p$& &$5$& &$7$& &$11$& &$13$& &$19$&  &$23$& &$31$& &$43$& &$67$& &$79$&
\cr \hh\om\om\om\om\om\om\om\om\om\om\om
\cr\nh\cr \hh\om\om\om\om\om\om\om\om\om\om\om
\cr &$g$& &$3$& &$5$& &$8$& &$11$& &$13$& &$21$& &$12$& &$34$& &$41$& &$53$&
\cr \hh\om\om\om\om\om\om\om\om\om\om\om
\cr\nh\cr \hh\om\om\om\om\om\om\om\om\om\om\om
\cr &$2^g-1\ \t{mod}\ p$& &$2$& &$3$& &$2$& &$6$& &$2$& &$11$& &$3$& &$20$& &$11$& &$30$&
\cr \hh\om\om\om\om\om\om\om\om\om\om\om
\cr\nh\cr \hh\om\om\om\om\om\om\om\om\om\om\om
\cr &$(g-1)!\ \t{mod}\ p$& &$2$& &$3$& &$2$& &$6$& &$10$& &$11$& &$22$& &$29$& &$44$& &$47$&
\cr \hh\om\om\om\om\om\om\om\om\om\om\om
\cr \nh\cr }}}
\smallskip

\begin{conjecture}\label{Conj3.4} {\rm (2014-05-11)} For any odd prime $p$, there exists a prime $q<p$ such that both $q$ and $2^q-q$
are primitive roots modulo $p$.
\end{conjecture}
\begin{remark}\label{Rem3.4}\rm We have verified this conjecture for odd primes $p<10^8$; see \cite[A242345]{S} for related data and graphs.
\end{remark}

\newpage
\centerline{Table 3.4: Primes $p$ with unique prime $q<p$ such that}
\centerline {both $q$ and $2^q-q$ are primitive roots modulo $p$}
\smallskip
\centerline{\vbox{\offinterlineskip
\halign{\vrule#&\ \ #\ \hfill   &&\vrule#&\ \ \hfill#\ \
\cr\nh\cr \hh\om\om\om
\cr &$p$& &$q<p$& &$2^q-q\ \t{mod}\ p$&
\cr \hh\om\om\om
\cr\nh\cr \hh\om\om\om
\cr &$3$& &$2$& &$2$&
\cr \hh\om\om\om
\cr\nh\cr \hh\om\om\om
\cr &$5$& &$2$& &$2$&
\cr \hh\om\om\om
\cr\nh\cr \hh\om\om\om
\cr &$7$& &$3$& &$5$&
\cr \hh\om\om\om
\cr\nh\cr \hh\om\om\om
\cr &$11$& &$2$& &$2$&
\cr \hh\om\om\om
\cr\nh\cr \hh\om\om\om
\cr &$13$& &$2$& &$2$&
\cr \hh\om\om\om
\cr\nh\cr \hh\om\om\om
\cr &$19$& &$2$& &$2$&
\cr \hh\om\om\om
\cr\nh\cr \hh\om\om\om
\cr &$23$& &$19$& &$7$&
\cr \hh\om\om\om
\cr\nh\cr \hh\om\om\om
\cr &$29$& &$2$& &$2$&
\cr \hh\om\om\om
\cr\nh\cr \hh\om\om\om
\cr &$31$& &$11$& &$22$&
\cr \hh\om\om\om
\cr\nh\cr \hh\om\om\om
\cr &$43$& &$3$& &$5$&
\cr \hh\om\om\om
\cr\nh\cr \hh\om\om\om
\cr &$61$& &$2$& &$2$&
\cr \hh\om\om\om
\cr\nh\cr \hh\om\om\om
\cr &$71$& &$67$& &$13$&
\cr \hh\om\om\om
\cr\nh\cr \hh\om\om\om
\cr &$73$& &$31$& &$58$&
\cr \hh\om\om\om
\cr\nh\cr \hh\om\om\om
\cr &$79$& &$59$& &$29$&
\cr \hh\om\om\om
\cr\nh\cr \hh\om\om\om
\cr &$97$& &$71$& &$74$&
\cr \hh\om\om\om
\cr\nh\cr \hh\om\om\om
\cr &$127$& &$43$& &$86$&
\cr \hh\om\om\om
\cr\nh\cr \hh\om\om\om
\cr &$151$& &$71$& &$14$&
\cr \hh\om\om\om
\cr \nh\cr }}}
\smallskip
\bigskip

Both Conjecture \ref{Conj3.4} and the following conjecture are more sophisticated than Erd\H os' Problem mentioned in Section 1.

\begin{conjecture}\label{qq!}
{\rm (i) (2014-04-21)} For any prime $p > 7$, there exists a prime $q < p$ such that both $q$ and $q!$ are primitive roots modulo $p$.

{\rm (ii) (2017-08-27)} For any odd prime $p$, there exists a prime $q<p$ such that $q$ is not only
a primitive root modulo $p$ but also a primitive root modulo $p_q$.
\end{conjecture}
\begin{remark}\label{Rem-qq!}\rm (a) We have verified this conjecture for primes $p<10^8$; see \cite[A236306 and A291615]{S} for related data and graphs. For example, both $3$ and $3!=6$
are primitive roots modulo the prime $17$; the number $3$ is a primitive root modulo the prime $43$ and also a primitive root modulo $p_3=5$.

 (b) If there are only finitely many primes $q$
with $q$ a primitive root modulo $p_q$, then for the product $P$ of all such primes $q$,
by Dirichlet's theorem on primes in arithmetic progressions, $p\eq1\pmod{4P}$ for some prime $p$,
hence for any prime $q\mid P$ we have $(\f qp)=1$ (by the law of quadratic reciprocity) and thus $q$ is not a primitive root modulo $p$. So Part (ii) of Conjecture \ref{Conj3.4} implies that there are infinitely primes $q$ such that $q$ is a primitive root modulo $p_q$. For such primes $q$, see \cite[A291657]{S}.
\end{remark}

\begin{conjecture}\label{Conj-FL} {\rm (i)} For any prime $p$, there are positive integers $k$ and $m$
such that $g=F_kF_m$ is smaller than $p$ and also a primitive root modulo $p$.
Moreover, the set $G(p)$ given by \eqref{G(p)} contains a number of the form $\ell F_m$ with $m\in\Z^+$ and
$\ell\in\{F_k:\ k=2,3,4\}=\{1,2,3\}$.

{\rm (ii)} For each prime $p$, there are $k,m\in\N$
such that $g=L_kL_m$ is smaller than $p$ and also a primitive root modulo $p$.
\end{conjecture}
\begin{remark}\label{FL} We have verified parts (i) and (ii) of Conjecture \ref{Conj-FL} for primes $p$ smaller than
$5\times10^9$ and $10^9$ respectively. See \cite[A331506]{S} for related data.
In contrast with Conjectures 2.1(i) and Conjecture 2.2(i), Conjecture \ref{Conj-FL} implies that for any prime $p$ there are $k,m\in\N$ with $F_k<p$ and $L_m<p$ such that
$(\f{F_k}p)=(\f{L_m}p)=-1$.
\end{remark}

\begin{conjecture}\label{5-10} {\rm (2018-05-24)} For any prime $p>7$, there is a a number $g=5^k+10^m$ with $k,m\in\N$
such that $g$ is smaller than $p$ and also a primitive root modulo $p$.
\end{conjecture}
\begin{remark} \label{5-10}
We have verified this for all primes $7<p<10^9$ (cf. \cite[A305048]{S}). Our computation suggests that
\begin{gather*}3,\ 5,\, 31,\ 43,\ 241,\ 307,\ 311,\ 421,\ 523,\ 547,\ 607,\ 727,\ 2311,\ 2511,
\\2689,\ 6091,\ 9439,\ 13381,\ 55441,\ 56401,\ 66301,\ 6276271
\end{gather*}
are the only values of primes $p$ which has a unique primitive root $g<p$ of the form $5^k+10^m\ (k,m\in\N)$. For example, for the prime $p=6276271$,  the unique primitive root $g<p$ of the form $5^k+10^m$ is $5^5+10=3135$.
\end{remark}

\begin{conjecture}\label{cen} {\rm (2018-05-24)}
{\rm (i)} For any odd prime $p$, the set $G(p)$ given by \eqref{G(p)} contains a number of the form $\bi{2k}k+\bi{2m}m$ with $k,m\in\N$.

{\rm (ii)} For any odd prime $p$, the set $G(p)$ given by \eqref{G(p)} contains a number of the form $C_k+C_m$ with $k,m\in\N$, where $C_n$ denotes the $n$-th Catalan number $\bi{2n}n/(n+1)$.
\end{conjecture}
\begin{remark} \label{cen}
We have verified Conjecture \ref{cen} for all odd primes $p<10^9$ (cf. \cite[A305030]{S}).
For example,
$\bi{2}1+\bi{8}{4}=2+70=72$
is a primitive root modulo the prime $109$.
\end{remark}

Recall that the Bernoulli numbers $B_0,B_1,B_2,\ldots$ are rational numbers defined by
$$B_0=1,\quad\t{and}\quad \sum_{k=0}^n\bi{n+1}kB_k=0\ \ \t{for all}\ n=1,2,3,\ldots,$$
and the Euler numbers $E_0,E_1,E_2,\ldots$ are integers defined by
$$E_0=1,\quad\t{and}\quad \sum^n_{k=0\atop2\mid n-k} \bi nkE_k=0\ \ \t{for all}\ n=1,2,3,\ldots.$$
It is well known that $B_{2n+1}=E_{2n-1}=0$ for all $n=1,2,3,\ldots$.
For any prime $p>3$ it is well known that all the Bernoulli numbers
$B_{2k}\ (k=1,\ldots,(p-3)/2)$
are $p$-adic integers (this follows from the recurrence for Bernoulli numbers
or Kummer's theorem on Bernoulli numbers). The tangent numbers $t_1,t_2,\ldots$ are given by
$$\tan x= \sum_{n=1}^\infty t_n\f{x^{2n-1}}{(2n-1)!}\ \ \text{for}\ |x|<\f{\pi}2. $$
It is known that
$$t_n=(-1)^{n-1}2^{2n}(2^{2n}-1)\f{B_{2n}}{2n}\qquad \ \t{for all}\ n\in\Z^+.$$

\begin{conjecture}\label{Conj-BET} {\rm (2014-05-07) (i)} For any prime $p>3$, there exists a prime $q<p$ such that the Bernoulli number $B_{q-1}$
is a primitive root modulo $p$.

{\rm (ii)} For any prime $p>13$, there exists a prime $q<p$ such that the Euler number $E_{q-1}$
is a primitive root modulo $p$.

{\rm (iii)} For any odd prime $p$, there is a prime $q<p$ such that the tangent number $t_q$
is a primitive root modulo $p$.
\end{conjecture}
\begin{remark}\label{Rem3.5} We have verified Conjecture \ref{Conj-BET} for primes $p<10^8$; see \cite[A242210 and A242213]{S} for related data and graphs.
\end{remark}

\medskip
\centerline{Table 3.5: Primes $p$ with unique prime $q<p$ such that}
\centerline{$B_{q-1}$ is a primitive root modulo $p$}
\smallskip
\centerline{\vbox{\offinterlineskip
\halign{\vrule#&\ \ #\ \hfill   &&\vrule#&\ \ \hfill#\ \
\cr\nh\cr \hh\om\om\om\om
\cr &$p$& &$5$& &$11$& &$19$&
\cr \hh\om\om\om\om
\cr\nh\cr \hh\om\om\om\om
\cr &$q<p$& &$2$& &$3$& &$17$&
\cr \hh\om\om\om\om
\cr\nh\cr \hh\om\om\om\om
\cr &$B_{q-1}$& &$-1/2$& &$1/6$& &$-3617/510$&
\cr \hh\om\om\om\om
\cr\nh\cr \hh\om\om\om\om
\cr &$B_{q-1}\ \t{mod}\ p$& &$2$& &$2$& &$15$&
\cr \hh\om\om\om\om
\cr \nh\cr }}}
\smallskip
\bigskip

Recall that those rational numbers $H_n=\sum_{0<k\ls n}1/k\ (n=0,1,2,\ldots)$ are called harmonic numbers.
The second-order harmonic numbers are those rational numbers $H_n^{(2)}=\sum_{0<k\ls n}1/k^2$ with $n\in\N$.

\begin{conjecture}\label{Conj-H} {\rm (2014-05-08)} Let $p>5$ be a prime.

{\rm (i)} There exists a prime $q\ls(p+1)/2$ such that $H_{q-1}$ is a primitive root modulo $p$.

{\rm (ii)} There exists a prime $q\ls(p-1)/2$ such that $H^{(2)}_{q-1}$ is a primitive root modulo $p$.
\end{conjecture}
\begin{remark}\label{Rem-H} We have verified both parts of Conjecture \ref{Conj-H} for all primes
$5<p<10^8$; See \cite[A242222 and A242241]{S} for related data and graphs.
\end{remark}

\begin{conjecture}\label{Conj-mix}
{\rm (i) (2014-04-21)} For any prime $p > 3$, there exists a prime $q < p/2$ such that the
Catalan number $C_q = \bi{2q}q/(q+1)$ is a primitive root modulo $p$.

{\rm  (ii) (2014-04-22)} For any prime $p > 3$, there exists a prime $q < p/2$ such that the
Bell number ${\rm Bell}(q)$ is a primitive root modulo $p$, where ${\rm Bell}(q)$ denotes the number of ways to partition a set of cardinality $q$.

{\rm (iii) (2014-05-11)} For any prime $p>3$, there exists a prime $q < p/2$ such that the
Franel number $f_{q-1}=\sum_{k=0}^q\bi {q-1}k^3$ is a primitive root modulo $p$.

{\rm (iv)} For any prime $p>3$, there exists a prime $q < p$ such that the
 $T_q$ is a primitive root modulo $p$, where the central trinomial coefficient $T_q$
 denotes the coefficient of $x^q$ in the expansion of $(x^2+x+1)^q$.

\end{conjecture}
\begin{remark}\label{Rem3.10}\rm  We have verified all parts of Conjecture \ref{Conj-mix}
for each prime $3<p<10^8$. For related data and graphs concerning parts (i)-(ii) of Conjecture
\ref{Conj-mix}, one may visit \cite[A236308 and A237594]{S}.
\end{remark}

\vskip 0.5cm
 \noindent{\Large\bf 4\quad Other Conjectures involving Primitive Roots Modulo Primes}
 \vskip 0.5cm

\renewcommand{\theequation}{4.\arabic{equation}}
\renewcommand{\thetheorem}{4.\arabic{theorem}}
\renewcommand{\thecorollary}{4.\arabic{corollary}}
\renewcommand{\thelemma}{4.\arabic{lemma}}
\renewcommand{\theconjecture}{4.\arabic{conjecture}}
\renewcommand{\theremark}{4.\arabic{remark}}

 \setcounter{lemma}{0} \setcounter{theorem}{0}
\setcounter{corollary}{0}
\setcounter{remark}{0}
\setcounter{equation}{0}
\setcounter{conjecture}{0}

The following conjecture was originally motivated by the Chinese Remainder Theorem.
\begin{conjecture}\label{Conj4.1} {\rm (2017-08-29)} {\rm (i)} Let $p$ and $q$ be primes.
Then there is a positive integer $g\ls\sqrt{4pq+1}$ such that $g$ is a primitive root modulo $p$ and also a primitive root modulo $q$. We may require further that $g <\sqrt{pq}$ unless $\{p,q\}$ is among the $15$ pairs
\begin{gather*}\{2,3\}, \ \{2,11\},\ \{2,13\},\ \{2,59\},\  \{2,131\},\ \{2,181\},
 \\\{3,7\},\ \{3,31\},\  \{3,79\},\ \{3,191\},\ \{3,199\},\ \{5,271\},\ \{7,11\},\ \{7,13\},\ \{7,71\}.
 \end{gather*}

{\rm (ii)} Let $n$ be any positive integer. If $q_1,\ldots,q_n$ are primes with $\max\{q_1,\ldots,q_n\}$ sufficiently large, then there is a positive integer $g \ls n!(q_1\ldots q_n)^{1/n}$  which is a primitive root modulo $q_k$ for all $k=1,\ldots,n$.
\end{conjecture}
\begin{remark}\label{Rem4.1}\rm See \cite[A291690]{S} for related data and comments. We have verified part (i) of Conjecture \ref{Conj4.1}
for primes $p,q<2\times10^5$. For example, $5=\sqrt{4\times2\times3+1}$ is the least positive integer which is a primitive root
modulo $2$ and also a primitive root modulo $3$, and $19=\lfloor\sqrt{4\times7\times13+1}\rfloor$
is the least positive integer which is a primitive root modulo $7$ and also a primitive root modulo $13$. Our computation for primes smaller than $3515$ suggests that if $q_1\ls q_2\ls q_3$ are primes but there is no positive integer
$g\ls 6\root 3\of{q_1q_2q_3}$ which is a primitive root modulo $q_i$ for all $i=1,2,3$, then
$(q_1,q_2,q_3)$ must be among the following 13 triples:
\begin{gather*}(3,5,43),\ (3,7,13),\ (3,7,19),\ (3,7,67),\ (3,7,127),\ (3,7,151),\ (3,7,421),
\\ (3,13,127),\ (3,31,43),\ (5,13,31),\ \ (7,11,523),\ (7,23,127),\ (31,37,79).
\end{gather*}
For $(q_1,q_2,q_3,q_4)=(3,31,43,991)$, $1439$ is the least positive integer which is a primitive root
modulo $q_j$ for all $j=1,2,3,4$. Note that $1439/(3\times31\times43\times991)^{1/4}\approx 32.25$.
\end{remark}

\begin{conjecture}\label{Conj4.2} {\rm (2015-08-05) (i)} For any prime $p>13$, there are distinct positive integers $a$ and $b$
with $a+b<p$ such that $a,b,a+b,ab(a+b)$ are all primitive roots modulo $p$.

{\rm (ii)} For any prime $p>13$ with $p\not=31$, there are $a,b,c\in\{1,\ldots,p-1\}$
with $a^2+b^2=c^2$ such that $abc$ is a primitive root modulo $p$.
\end{conjecture}
\begin{remark}\label{Rem4.2} We have verified parts (i) and (ii) for primes below $10^6$ and $10^8$ respectively. See \cite[A260947 and A260946]{S} for related data.
\end{remark}

\begin{conjecture}\label{Conj-tri} {\rm (2015-08-06)} {\rm (i)} For any prime $p > 7$, there exists a right triangle whose three sides are among $1,\ldots,p-1$ and whose area is a primitive root modulo $p$.

{\rm (ii)} For any prime $p > 31$, there exists a right triangle whose three sides are among $1,\ldots,p-1$, and whose perimeter and area are quadratic residues modulo $p$.
\end{conjecture}
\begin{remark} See \cite[A260960]{S} for related data and graphs. For example, $6$ is a primitive
root modulo the prime $17$ and $6$ is also the area of a right triangle with sides $3,4,5$.
\end{remark}

\begin{conjecture}\label{Conj4.3} {\rm (2014-06-11)} For any $m,N\in\Z^+$, there is a positive integer $n\gs N$ such that $p_{n+i}$ is a primitive root modulo $p_{n+j}$ for all $i,j=0,\ldots,m$ with $i\not=j$.
\end{conjecture}
\begin{remark} See \cite[A243839]{S} for related data.
Via the Maynard-Tao theorem, H. Pan and Z.-W. Sun \cite{PS} showed that the Generalized Riemann Hypothesis implies Conjecture \ref{Conj4.3}.
\end{remark}

\begin{conjecture}\label{Conj-phi} {\rm (2017-10-02)} There are infinitely many primes $p$ such that
$\phi(p-1)$ is a primitive root modulo $p$. Moreover, there is a constant $0.361<s<0.362$ such that
$\lim_{x\to+\infty}{S(x)}/(x/\log x)=s,$
where $S(x)$ denotes the number of primes $p\ls x$ with $\varphi(p-1)$ a primitive root modulo $p$.
\end{conjecture}
\begin{remark} It is well known that for any prime $p$ there are exactly $\varphi(p-1)$ numbers among $1,\ldots,p-1$ which are primitive roots modulo $p$.
See \cite[A293213]{S} for such special primes $p$ with $\varphi(p-1)$ a primitive root modulo $p$.
Among the first $6\times10^8$ primes, there are exactly $216635723$ such special primes. Note that
$216635723/(6\times10^8)\approx 0.36105954$.
\end{remark}

\begin{conjecture}\label{Conj-three-g} {\rm (2013-10-02)} For any prime $p > 7$ with $p\not=13,29,61$, there are three consecutive integers among $1,\ldots,p-1$ which are primitive roots modulo $p$.
\end{conjecture}
\begin{remark}
See \cite[A229899]{S} for related data. For example, $19,20,21$ are primitive roots modulo the prime $23$.
\end{remark}

\begin{conjecture}\label{Conj-gap} {\rm (2017-10-01)} Let $p$ be a prime. For any $x\in\Z$, one of $x,x+1,\ldots,x+2\lfloor\sqrt{p+2}\rfloor+2$ is a primitive root modulo $p$.
\end{conjecture}
\begin{remark} We have verified this for all primes $p<10^5$. For example, $2\lfloor\sqrt{79+2}\rfloor+2=20$
and among the integers $8+k\ (k=0,\ldots,20)$ only $28$ is a primitive root modulo $79$.
Also, $2\lfloor\sqrt{409+2}\rfloor+2=42$, and among the integers $388+k\ (k=0,\ldots,42)$
only $388$ and $430$ are primitive roots modulo the prime $409$.
By \cite[Theorem 3]{Bu}, Conjecture \ref{Conj-gap} holds for sufficiently large primes $p$.
\end{remark}

For $a,b,c\in\Z$, we set
 \begin{equation}\label{Sabc}S_p(a,b,c):=\sum_{g\in G(p)}\l(\f{ag^2+bg+c}p\r),
 \end{equation}
 where $G(p)$ is given by \eqref{G(p)}.
  Since the inverse $g^*$ of $g\in G(p)$ modulo $p$ is also a primitive root modulo $p$,
 we see that
 $$S_p(a,b,c)=\sum_{g\in G(p)}\l(\f{a(g^*)^{2}+bg^*+c}p\r)=\sum_{g\in G(p)}\l(\f{a+bg+cg^2}p\r)=S_p(c,b,a).$$

\begin{conjecture}\label{Conj-g} {\rm (2013-10-02)} Let $p>11$ be a prime, and let $a,b,c\in\Z$ with $b^2-4ac\not\eq0\pmod p$. If $a$ or $c$ is not divisible by $p$, then
\begin{equation}\label{bound}|S_p(a,b,c)|<\f{\sqrt p}2\log p.
\end{equation}
\end{conjecture}
\begin{remark} We note that $S_{11}(1,-3,1)/(\sqrt{11}\log 11)\approx 0.50296$.
\end{remark}

\begin{conjecture}\label{g-quad} {\rm (2013-10-02)} Let $p > 13$ be a prime with $p\not=19,31$, and let $a,b,c$ be integers with $a$ or $c$ not divisible by $p$. If $p$ does not divide $b^2-4ac$, then there is a primitive root $g$ modulo $p$ such that $ag^2+bg+c$ is a quadratic residue modulo $p$, and there is also a primitive root $h$ modulo $p$ such that $ah^2+bh+c$ is a quadratic nonresidue modulo $p$.
\end{conjecture}
\begin{remark}\label{g-quad} Compare this with a consequence of Theorem 1.1 mentioned in Section 1.
\end{remark}

\begin{conjecture}\label{g-1} {\rm (2014-04-20)} Let $a$ be any positive integer.

{\rm (i)} For any prime $p>7$ with $p\not =13$, there is a primitive root $g$ modulo $p^a$ such that
$g+g^{-1}$ is also a primitive root modulo $p^a$.

{\rm (ii)} For any prime $p>5$ with $p\not =13,61$, there is a primitive root $g$ modulo $p^a$ such that
$g-g^{-1}$ is also a primitive root modulo $p^a$.
\end{conjecture}
\begin{remark}\label{Rem-g-1} Note that Conjecture \ref{g-1} holds for $a=1$ by Theorem \ref{Th1.2}.
\end{remark}

When $p$ and $2p+1$ are both prime, $p$ is called a Sophie Germain prime.

\begin{conjecture}\label{Conj4.4} {\rm (2014-01-17)} {\rm (i)} Any integer $n>37$ can be written as
$k+m$ with $k,m\in\Z^+$ such that $p=p_k+\varphi(m)$ is a Sophie Germain prime having $2$ as a primitive root.

{\rm (ii)} Any integer $n>7$ can be written as
$k+m$ with $k,m\in\Z^+$ such that $p=\varphi(k)+\varphi(m/2)-1$ is a prime having $2$ as a primitive root.
\end{conjecture}
\begin{remark} See \cite[A235987]{S} for related data and graphs.
\end{remark}

\vskip 0.5cm
 \noindent{\Large\bf 5\quad Primitive Prime Divisors of Some Combinatorial Sequences}
 \vskip 0.5cm

\renewcommand{\theequation}{5.\arabic{equation}}
\renewcommand{\thetheorem}{5.\arabic{theorem}}
\renewcommand{\thecorollary}{5.\arabic{corollary}}
\renewcommand{\thelemma}{5.\arabic{lemma}}
\renewcommand{\theconjecture}{5.\arabic{conjecture}}
\renewcommand{\theremark}{5.\arabic{remark}}

 \setcounter{lemma}{0} \setcounter{theorem}{0}
\setcounter{corollary}{0}
\setcounter{remark}{0}
\setcounter{equation}{0}
\setcounter{conjecture}{0}

\begin{conjecture}\label{Conj5.1} For any integer $n>1$ with $n\not=5,16$, the number $2^n-n$ has a prime divisor $p$
not dividing any $2^k-k$ with $0<k<n$.
\end{conjecture}
\begin{remark}\label{Rem5.1}\rm See \cite[A242292]{S} for related data.
\end{remark}

\begin{conjecture}\label{Conj5.2}  For any integer $n>4$, there is a prime $p$ for which $B_{2n}\eq0\pmod p$ but
$B_{2k}\not\eq0\pmod p$ for all $0<k<n$. For each $n=2,3,\ldots$, the Euler number $E_{2n}$
has a prime divisor $p$ not dividing any $E_{2k}$ with $0<k<n$. Also, for every $n=4,5,\ldots$
the tangent number $t_n$ has a prime divisor $p$ not dividing any $t_k$ with $0<k<n$.
\end{conjecture}
\begin{remark}\label{Rem5.2}\rm See \cite[A242193, A242194 and A242195]{S}.
In Table 5.1, $p_B(n)$ denotes the least prime $p$ for which $B_{2n}\eq0\pmod p$ but
$B_{2k}\not\eq0\pmod p$ for all $0<k<n$, similarly $p_E(n)$ represents the least prime divisor of $E_{2n}$ not dividing any $E_{2k}$ with $0<k<n$.
\end{remark}

\newpage
\centerline{Table 5.1: Least primitive prime divisors $p_B(n)$ of $B_{2n}$ and $p_E(n)$ of $E_{2n}$}
\smallskip
\centerline{\vbox{\offinterlineskip
\halign{\vrule#&\ \ #\ \hfill   &&\vrule#&\ \ \hfill#\ \
\cr\nh\cr \hh\om\om\om
\cr &$n$& &$p_B(n)$& &$p_E(n)$&
\cr \hh\om\om\om
\cr\nh\cr \hh\om\om\om
\cr &$2$& &$ $& &$5$&
\cr \hh\om\om\om
\cr\nh\cr \hh\om\om\om
\cr &$3$& &$ $& &$61$&
\cr \hh\om\om\om
\cr\nh\cr \hh\om\om\om
\cr &$4$& &$ $& &$277$&
\cr \hh\om\om\om
\cr\nh\cr \hh\om\om\om
\cr &$5$& &$5$& &$19$&
\cr \hh\om\om\om
\cr\nh\cr \hh\om\om\om
\cr &$6$& &$691$& &$13$&
\cr \hh\om\om\om
\cr\nh\cr \hh\om\om\om
\cr &$7$& &$7$& &$47$&
\cr \hh\om\om\om
\cr\nh\cr \hh\om\om\om
\cr &$8$& &$3617$& &$17$&
\cr \hh\om\om\om
\cr\nh\cr \hh\om\om\om
\cr &$9$& &$43867$& &$79$&
\cr \hh\om\om\om
\cr\nh\cr \hh\om\om\om
\cr &$10$& &$283$& &$41737$&
\cr \hh\om\om\om
\cr\nh\cr \hh\om\om\om
\cr &$11$& &$11$& &$31$&
\cr \hh\om\om\om
\cr\nh\cr \hh\om\om\om
\cr &$12$& &$103$& &$2137$&
\cr \hh\om\om\om
\cr\nh\cr \hh\om\om\om
\cr &$13$& &$13$& &$67$&
\cr \hh\om\om\om
\cr\nh\cr \hh\om\om\om
\cr &$14$& &$9349$& &$29$&
\cr \hh\om\om\om
\cr\nh\cr \hh\om\om\om
\cr &$15$& &$1721$& &$15669721$&
\cr \hh\om\om\om
\cr\nh\cr \hh\om\om\om
\cr &$16$& &$37 $& &$930157$&
\cr \hh\om\om\om
\cr\nh\cr \hh\om\om\om
\cr &$17$& &$17$& &$4153$&
\cr \hh\om\om\om
\cr\nh\cr \hh\om\om\om
\cr &$18$& &$26315271553053477373$& &$37$&
\cr \hh\om\om\om
\cr\nh\cr \hh\om\om\om
\cr &$19$& &$19$& &$23489580527043108252017828576198947741$&
\cr \hh\om\om\om
\cr\nh\cr \hh\om\om\om
\cr &$20$& &$137616929$& &$41$&
\cr \hh\om\om\om
\cr\nh\cr \hh\om\om\om
\cr &$21$& &$1520097643918070802691$& &$137$&
\cr \hh\om\om\om
\cr\nh\cr \hh\om\om\om
\cr &$22$& &$59$& &$587$&
\cr \hh\om\om\om
\cr\nh\cr \hh\om\om\om
\cr &$23$& &$23$& &$285528427091$&
\cr \hh\om\om\om
\cr\nh\cr \hh\om\om\om
\cr &$24$& &$653$& &$5516994249383296071214195242422482492286460673697$&
\cr \hh\om\om\om
\cr\nh\cr \hh\om\om\om
\cr &$25$& &$417202699$& &$5639$&
\cr \hh\om\om\om
\cr\nh\cr \hh\om\om\om
\cr &$26$& &$577$& &$53$&
\cr \hh\om\om\om
\cr\nh\cr \hh\om\om\om
\cr &$27$& &$39409$& &$2749$&
\cr \hh\om\om\om
\cr\nh\cr \hh\om\om\om
\cr &$28$& &$113161$& &$5303$&
\cr \hh\om\om\om
\cr\nh\cr \hh\om\om\om
\cr &$29$& &$29$& &$1459879476771247347961031445001033$&
\cr \hh\om\om\om
\cr\nh\cr \hh\om\om\om
\cr &$30$& &$2003$& &$6821509$&
\cr \hh\om\om\om
\cr\nh\cr \hh\om\om\om
\cr &$31$& &$31$& &$101$&
\cr \hh\om\om\om
\cr\nh\cr \hh\om\om\om
\cr &$32$& &$1226592271$& &$25349$&
\cr \hh\om\om\om
\cr \nh\cr }}}
\smallskip
\bigskip

\begin{conjecture}\label{Conj5.3} Let $n>1$ be an integer. If $n\not=7$, then there is a prime $p$ for which $H_n\eq0\pmod p$
but $H_k\not\eq0\pmod p$ for all $0<k<n$.
Also, there is a prime $p$ for which $H_n^{(2)}\eq0\pmod p$
but $H_k^{(2)}\not\eq0\pmod p$ for all $0<k<n$.
\end{conjecture}
\begin{remark}\label{Rem5.3}\rm For related numerical data, see \cite[A242223 and A242241]{S}.
\end{remark}

\centerline{Table 5.2: Least primitive prime divisors $p_H(n)$ of $H_{n}$ and $p_{H^{(2)}}(n)$ of $H^{(2)}_n$}
\smallskip
\centerline{\vbox{\offinterlineskip
\halign{\vrule#&\ \ #\ \hfill   &&\vrule#&\ \ \hfill#\ \
\cr\nh\cr \hh\om\om\om
\cr &$n$& &$p_H(n)$& &$p_{H^{(2)}}(n)$&
\cr \hh\om\om\om
\cr\nh\cr \hh\om\om\om
\cr &$2$& &$3$& &$5$&
\cr \hh\om\om\om
\cr\nh\cr \hh\om\om\om
\cr &$3$& &$11$& &$7$&
\cr \hh\om\om\om
\cr\nh\cr \hh\om\om\om
\cr &$4$& &$5$& &$41$&
\cr \hh\om\om\om
\cr\nh\cr \hh\om\om\om
\cr &$5$& &$137$& &$11$&
\cr \hh\om\om\om
\cr\nh\cr \hh\om\om\om
\cr &$6$& &$7$& &$13$&
\cr \hh\om\om\om
\cr\nh\cr \hh\om\om\om
\cr &$7$& &$ $& &$266681$&
\cr \hh\om\om\om
\cr\nh\cr \hh\om\om\om
\cr &$8$& &$761$& &$17$&
\cr \hh\om\om\om
\cr\nh\cr \hh\om\om\om
\cr &$9$& &$7129$& &$19$&
\cr \hh\om\om\om
\cr\nh\cr \hh\om\om\om
\cr &$10$& &$61$& &$178939$&
\cr \hh\om\om\om
\cr\nh\cr \hh\om\om\om
\cr &$11$& &$97$& &$23$&
\cr \hh\om\om\om
\cr\nh\cr \hh\om\om\om
\cr &$12$& &$13$& &$18500393$&
\cr \hh\om\om\om
\cr\nh\cr \hh\om\om\om
\cr &$13$& &$29$& &$40799043101$&
\cr \hh\om\om\om
\cr\nh\cr \hh\om\om\om
\cr &$14$& &$1049$& &$29$&
\cr \hh\om\om\om
\cr\nh\cr \hh\om\om\om
\cr &$15$& &$41233$& &$31$&
\cr \hh\om\om\om
\cr\nh\cr \hh\om\om\om
\cr &$16$& &$17 $& &$619$&
\cr \hh\om\om\om
\cr\nh\cr \hh\om\om\om
\cr &$17$& &$37$& &$601$&
\cr \hh\om\om\om
\cr\nh\cr \hh\om\om\om
\cr &$18$& &$19$& &$8821$&
\cr \hh\om\om\om
\cr\nh\cr \hh\om\om\om
\cr &$19$& &$7440427$& &$86364397717734821$&
\cr \hh\om\om\om
\cr\nh\cr \hh\om\om\om
\cr &$20$& &$11167027$& &$421950627598601$&
\cr \hh\om\om\om
\cr\nh\cr \hh\om\om\om
\cr &$21$& &$18858053$& &$2621$&
\cr \hh\om\om\om
\cr\nh\cr \hh\om\om\om
\cr &$22$& &$23$& &$295831$&
\cr \hh\om\om\om
\cr\nh\cr \hh\om\om\om
\cr &$23$& &$583859$& &$47$&
\cr \hh\om\om\om
\cr\nh\cr \hh\om\om\om
\cr &$24$& &$577$& &$2237$&
\cr \hh\om\om\om
\cr\nh\cr \hh\om\om\om
\cr &$25$& &$109$& &$157$&
\cr \hh\om\om\om
\cr\nh\cr \hh\om\om\om
\cr &$26$& &$34395742267$& &$53$&
\cr \hh\om\om\om
\cr\nh\cr \hh\om\om\om
\cr &$27$& &$521$& &$307$&
\cr \hh\om\om\om
\cr\nh\cr \hh\om\om\om
\cr &$28$& &$375035183$& &$7741$&
\cr \hh\om\om\om
\cr\nh\cr \hh\om\om\om
\cr &$29$& &$4990290163$& &$6823$&
\cr \hh\om\om\om
\cr\nh\cr \hh\om\om\om
\cr &$30$& &$31$& &$61$&
\cr \hh\om\om\om
\cr\nh\cr \hh\om\om\om
\cr &$31$& &$2667653736673$& &$205883$&
\cr \hh\om\om\om
\cr\nh\cr \hh\om\om\om
\cr &$32$& &$2917$& &$487$&
\cr \hh\om\om\om
\cr \nh\cr }}}
\smallskip

In Table 5.2, $p_H(n)$ denotes the least prime $p$ for which $H_n\eq0\pmod p$
but $H_k\not\eq0\pmod p$ for all $0<k<n$, and $p_{H^{(2)}}(n)$ represents the least prime $p$ for which $H_n^{(2)}\eq0\pmod p$
but $H_k^{(2)}\not\eq0\pmod p$ for all $0<k<n$.

\begin{conjecture}\label{Conj-Be} For the sequence $({\rm Bell}(n))_{n>1}$ of Bell numbers,
each term ${\rm Bell}(n)$ with $n\gs2$ has a primitive prime divisor.
\end{conjecture}
\begin{remark} See \cite[A242171]{S} for related data.
\end{remark}

\centerline{Table 5.3: Least primitive prime divisors $p_{b}(n)$ of ${\rm Bell}(n)$ and $p_f(n)$ of $f_n=\sum_{k=0}^n\bi nk^3$}
\smallskip
\centerline{\vbox{\offinterlineskip
\halign{\vrule#&\ \ #\ \hfill   &&\vrule#&\ \ \hfill#\ \
\cr\nh\cr \hh\om\om\om
\cr &$n$& &$p_{b}(n)$& &$p_{f}(n)$&
\cr \hh\om\om\om
\cr\nh\cr \hh\om\om\om
\cr &$1$& &$ $& &$2$&
\cr \hh\om\om\om
\cr\nh\cr \hh\om\om\om
\cr &$2$& &$2$& &$5$&
\cr \hh\om\om\om
\cr\nh\cr \hh\om\om\om
\cr &$3$& &$5$& &$7$&
\cr \hh\om\om\om
\cr\nh\cr \hh\om\om\om
\cr &$4$& &$3$& &$173$&
\cr \hh\om\om\om
\cr\nh\cr \hh\om\om\om
\cr &$5$& &$13$& &$563$&
\cr \hh\om\om\om
\cr\nh\cr \hh\om\om\om
\cr &$6$& &$7$& &$13$&
\cr \hh\om\om\om
\cr\nh\cr \hh\om\om\om
\cr &$7$& &$877$& &$41$&
\cr \hh\om\om\om
\cr\nh\cr \hh\om\om\om
\cr &$8$& &$23$& &$369581$&
\cr \hh\om\om\om
\cr\nh\cr \hh\om\om\om
\cr &$9$& &$19$& &$937$&
\cr \hh\om\om\om
\cr\nh\cr \hh\om\om\om
\cr &$10$& &$4639$& &$61$&
\cr \hh\om\om\om
\cr\nh\cr \hh\om\om\om
\cr &$11$& &$22619$& &$23$&
\cr \hh\om\om\om
\cr\nh\cr \hh\om\om\om
\cr &$12$& &$37$& &$29$&
\cr \hh\om\om\om
\cr\nh\cr \hh\om\om\om
\cr &$13$& &$27644437$& &$2141$&
\cr \hh\om\om\om
\cr\nh\cr \hh\om\om\om
\cr &$14$& &$1800937$& &$12148537$&
\cr \hh\om\om\om
\cr\nh\cr \hh\om\om\om
\cr &$15$& &$251 $& &$31$&
\cr \hh\om\om\om
\cr\nh\cr \hh\om\om\om
\cr &$16$& &$241$& &$157$&
\cr \hh\om\om\om
\cr\nh\cr \hh\om\om\om
\cr &$17$& &$255755771$& &$59$&
\cr \hh\om\om\om
\cr\nh\cr \hh\om\om\om
\cr &$18$& &$19463$& &$37$&
\cr \hh\om\om\om
\cr\nh\cr \hh\om\om\om
\cr &$19$& &$271$& &$506251$&
\cr \hh\om\om\om
\cr\nh\cr \hh\om\om\om
\cr &$20$& &$61$& &$151$&
\cr \hh\om\om\om
\cr\nh\cr \hh\om\om\om
\cr &$21$& &$24709$& &$3019$&
\cr \hh\om\om\om
\cr\nh\cr \hh\om\om\om
\cr &$22$& &$17$& &$769$&
\cr \hh\om\om\om
\cr\nh\cr \hh\om\om\om
\cr &$23$& &$89$& &$47$&
\cr \hh\om\om\om
\cr\nh\cr \hh\om\om\om
\cr &$24$& &$123419$& &$6730949$&
\cr \hh\om\om\om
\cr\nh\cr \hh\om\om\om
\cr &$25$& &$367$& &$79$&
\cr \hh\om\om\om
\cr\nh\cr \hh\om\om\om
\cr &$26$& &$101$& &$53$&
\cr \hh\om\om\om
\cr\nh\cr \hh\om\om\om
\cr &$27$& &$157$& &$3853$&
\cr \hh\om\om\om
\cr\nh\cr \hh\om\om\om
\cr &$28$& &$67$& &$661$&
\cr \hh\om\om\om
\cr\nh\cr \hh\om\om\om
\cr &$29$& &$75979$& &$138961158000728258971$&
\cr \hh\om\om\om
\cr\nh\cr \hh\om\om\om
\cr &$30$& &$107$& &$1361$&
\cr \hh\om\om\om
\cr \nh\cr }}}
\smallskip
\bigskip

Table 5.3 is related to Conjecture \ref{Conj-Be} and
the following Conjecture \ref{Conj-Fr}, where $p_{b}(n)$ denotes the least prime divisor $p$ of the Bell number ${\rm Bell}(n)$ which does not divide any ${\rm Bell}(k)$ with $0<k<n$,
and $p_{f}(n)$ represents the least prime divisor $p$ of the Franel number $f_n=\sum_{k=0}^n\bi nk^3$ which does not any $f_k$ with $0<k<n$.

\begin{conjecture}\label{Conj-Fr}
 For the sequence $(f_n)_{n\gs1}$ of Franel numbers, each term $f_n=\sum_{k=0}^n\bi nk^3$ with $n\in\Z^+$ has a primitive prime divisor.
For the sequence $(f_n^{(4)})_{n\gs1}$ of the fourth-order Franel numbers with $f_n^{(4)}=\sum_{k=0}^n\bi nk^4$, each term $f_n^{(4)}$ with $n\in\Z^+$ has a primitive prime divisor.
In general, for any integer $r>2$,  if $n\in\Z^+$ is large enough then $f_n^{(r)}=\sum_{k=0}^n\bi nk^r$ has a prime divisor $p$
not dividing any $f_k^{(r)}$ with $0<k<n$.
\end{conjecture}
\begin{remark}\label{Rem5.4}\rm For related numerical data, see \cite[A242171 and A242169]{S}.
\end{remark}

For each $n\in\N$, the central trinomial coefficient $T_n=\sum_{k=0}^{\lfloor n/2\rfloor}\bi n{2k}\bi{2k}k$
is the coefficient of $x^n$ in the expansion of $(x^2+x+1)^n$, and the Motzkin number $M_n$
is given by $M_n=\sum_{k=0}^n\bi n{2k}C_k$.

\begin{conjecture}\label{Conj-T}  {\rm (i)} For the sequence $(T_n)_{n\gs1}$ of central trinomial coefficients, each term $T_n$ with $n>1$ has a primitive prime divisor.

{\rm (ii)} Each term of the sequence $(M_n)_{n\gs4}$ of Motzkin numbers has a primitive prime divisor.
\end{conjecture}
\begin{remark}\rm See \cite[A242170]{S} for related data. For integer $n>1$ let $p_T(n)$ be the least prime factor of $T_n$ which does not divide any of $T_k\ (0<k<n)$. Then
\begin{gather*}p_T(2)=3,\ p_T(3)=7,\ p_T(4)=19,\ p_T(5)=17,\ p_T(6)=47,\ p_T(7)=131,
\ p_T(8)=41,
\\ p_T(9)=43,\ p_T(10)=1279,\ p_T(11)=503,\ p_T(12)=113,\ p_T(13)=2917,\ p_T(14)=569,
\\p_T(15)=198623,\, p_T(16)=14083,\, p_T(17)=26693,\, p_T(18)=201611,\, p_T(19)=42998951,
\\ p_T(20)=41931041,\ p_T(21)=52635749,\ p_T(22)=1296973,\  p_T(23)=169097,
\\p_T(24)=1451,\ p_T(25)=1304394227,\ p_T(26)=107,\ p_T(27)=233,\ p_T(28)=173.
\end{gather*}
\end{remark}

Recall that the central Delannoy numbers $D_n\ (n\in\N)$ and the Ap\'ery numbers $A_n\ (n\in\N)$
are given by
$$D_n:=\sum_{k=0}^n\bi nk\bi{n+k}k\ \ \t{and}\ \ A_n:=\sum_{k=0}^n\bi nk^2\bi{n+k}k^2.$$

\begin{conjecture} Each term of the sequence $(D_n)_{n\gs1}$ of central Delannoy numbers has a primitive
prime divisor. Also, any term of the sequence $(A_n)_{n\gs1}$ of Ap\'ery numbers has a primitive prime divisor.
\end{conjecture}
\begin{remark} See \cite[A242173]{S} for related data.
\end{remark}

\begin{conjecture} Each term of the sequence $(d_n)_{n\gs3}$ of derangement numbers has a primitive prime divisor, where $d_n:=n!\sum_{k=0}^n\f{(-1)^k}{k!}$. Also,
any term of the sequence $({\rm Domb}(n))_{n\gs4}$ of Domb numbers has a primitive
prime divisor, where ${\rm Domb}(n):=\sum_{k=0}^n\bi nk^2\bi{2k}k\bi{2(n-k)}{n-k}$.
\end{conjecture}
\begin{remark} See \cite[A242207]{S} for related data.
\end{remark}

\begin{conjecture} For $n\in\Z^+$ let $q(n)$ denote the number of unordered ways to write
$n$ as a sum of distinct positive integers. Then, for any integer $n>203$, the number $q(n)$
has a prime divisor $p$ not dividing any $q(k)$ with $0<k<n$.
\end{conjecture}
\begin{remark} See \cite[A242180]{S} for related data. It is known that $q(n)\sim e^{\pi\sqrt{n/3}}/(4\root 4\of{3n^3})$ as $n\to+\infty$.
\end{remark}

Finally, we mention that Conjectures 5.1-5.9 were formulated by the author on May 7, 2014 on the basis of related computations.

\Ack. The initial version of this paper was posted to arXiv in May 2014
 with the ID {\tt arXiv:1405.0290}.
The author would like to thank Prof. Carl Pomerance and J. Wu for helpful comments.

\end{document}